\documentclass{article}
\usepackage{amssymb,amsmath,amscd,latexsym}
\usepackage[all]{xy}           %xypic macro for latex2.09

\newcommand{\efootnote}[1]{}

\newcommand{\delete}[1]{}
\newcommand{\mlabel}[1]{\label{#1}}  % Use this line to suppress names

\def \addgap{\vspace{6pt}}    %% automatic spacing before theorems
\def \delgap{\vspace{-4pt}}   %% eliminate spacing before items
\def \addmed{\addgap}
\def \delmed{\delgap}
\def \addsm{\vspace{3pt}}     %% fine tuning (manual)
\def \delsm{\vspace{-3pt}}    %% fine tuning (manual)
\def \delbig{\vspace{-12pt}}  %% really bad vertical underfull tuning

\delete{
\def \addgap{\vspace{11pt}}
\def \addsm{\vspace{10pt}}
\def \delbig{\vspace{-13pt}}
\def \delgap{\vspace{-8pt}}
\def \delsm{\vspace{-4pt}}
}

\makeatletter
\def\@cite#1#2{[{#1\if@tempswa , #2\fi}]}
\makeatother

\makeatletter
\def\th@plain{\thm@preskip\smallskipamount\thm@postskip\smallskipamount\itshape}
\makeatother

\numberwithin{equation}{section}

\newtheorem{tempthm}[equation]{Theorem}
\newtheorem{tempdefn}[equation]{Definiton}
\newtheorem{tempprop}[equation]{Proposition}
\newtheorem{tempcor}[equation]{Corollary}
\newtheorem{templem}[equation]{Lemma}
\newtheorem{tempex}[equation]{Example}
\newtheorem{tempexs}[equation]{Examples}
\newtheorem{temprmk}[equation]{Remark}
\newtheorem{tempexer}[equation]{Exercise}
\newenvironment{thm}{\addmed\begin{tempthm}\sl}{\end{tempthm}}
\newenvironment{defn}{\addmed\begin{tempdefn}\rm}{\end{tempdefn}}
\newenvironment{prop}{\addmed\begin{tempprop}\sl}{\end{tempprop}}

\newenvironment{exam}{\addmed\begin{tempex}\rm}{\end{tempex}}
\newenvironment{exams}{\addmed\begin{tempexs}\rm}{\end{tempexs}}

\newenvironment{proof}{\addmed \noindent {\em Proof.\,}}{\hspace{\fill} $\square$ \addmed}

\def\been{\begin{enumerate}}
\def\enen{\end{enumerate}}
\def\bea{\begin{eqnarray}}
\def\eea{\end{eqnarray}}
\def\bes{\begin{eqnarray*}}
\def\ees{\end{eqnarray*}}
\def\bei{\begin{itemize}}
\def\eni{\end{itemize}}
\def\beq{\begin{equation}}
\def\eeq{\end{equation}}

\newcommand{\nc}{\newcommand}
\nc{\bincc}[2]{  \left ( {\scs{#1} \atop 
    \vspace{-1cm}\scs{#2}} \right )}  %binomial coeff  
\nc{\binc}{\bincc}
\nc{\bs}{\overline{S}}
\nc{\la}{\longrightarrow}
\nc{\rar}{\rightarrow}
\nc{\dar}{\downarrow}
\nc{\dap}[1]{\downarrow \rlap{$\scriptstyle{#1}$}}
\nc{\uap}[1]{\uparrow \rlap{$\scriptstyle{#1}$}}
\nc{\defeq}{:=} %{\stackrel{\rm def}{=}}
\nc{\dis}[1]{\displaystyle{#1}}
\nc{\dotcup}{\ \displaystyle{\bigcup^\bullet}\ }
\nc{\hcm}{\ \hat{,}\ }
\nc{\hts}{\widehat{\otimes}}
\nc{\hcirc}{\hat{\circ}}
\nc{\lleft}{[}
\nc{\lright}{]}
\nc{\curlyl}{\left \{ \begin{array}{c} {} \\ {} \end{array}
    \right .  \!\!\!\!\!\!\!} 
\nc{\curlyr}{ \!\!\!\!\!\!\!
    \left . \begin{array}{c} {} \\ {} \end{array}
    \right \} }
\nc{\longmid}{\left | \begin{array}{c} {} \\ {} \end{array}
    \right . \!\!\!\!\!\!\!}
\nc{\ora}[1]{\stackrel{#1}{\rar}}
\nc{\ola}[1]{\stackrel{#1}{\la}}%${\Bbb Z}$
\nc{\scs}[1]{\scriptstyle{#1}}
\nc{\sss}{\subsubsection}
\nc{\mrm}[1]{{\rm #1}}
\nc{\margin}[1]{\marginpar{\rm #1}}   %{\rm #1}}
\nc{\dirlim}{\displaystyle{\lim_{\longrightarrow}}\,}
\nc{\invlim}{\displaystyle{\lim_{\longleftarrow}}\,}
\nc{\mvp}{\vspace{0.3cm}}
\nc{\tk}{^{(k)}}
\nc{\tp}{^\prime}
\nc{\ttp}{^{\prime\prime}}
\nc{\proofbegin}{\noindent{\bf Proof: }}
\nc{\proofend}{$\blacksquare$ \vspace{0.3cm}}
\nc{\modg}[1]{\!<\!\!{#1}\!\!>}
\nc{\intg}[1]{F_C(#1)}
\nc{\lmodg}{\!<\!\!}
\nc{\rmodg}{\!\!>\!}
\nc{\cpi}{\widehat{\Pi}}
\nc{\sha}{\mbox{\cyr X}}  %used to be \cyr
\nc{\shpr}{\diamond}    %Shuffle product
\nc{\shprc}{\shpr_c}
\nc{\labs}{\mid\!}
\nc{\rabs}{\!\mid}
\nc{\lpair}[2]{{<{#1} \mid {#2}>_\lambda}}
\nc{\mnote}[1]{{Note: {#1}}}
\nc{\Cq}{{C\!\!<\!\! {\mathfrak q}\!\!>}}

\nc{\ann}{\mrm{ann}}
\nc{\Aut}{\mrm{Aut}}
\nc{\can}{\mrm{can}}
\nc{\colim}{\mrm{colim}}
\nc{\Cont}{\mrm{Cont}}
\nc{\rchar}{\mrm{char}}
\nc{\cok}{\mrm{coker}}
\nc{\dtf}{{R-{\rm tf}}}
\nc{\dtor}{{R-{\rm tor}}}

\nc{\Div}{{\mrm Div}}
\nc{\End}{\mrm{End}}
\nc{\Ext}{\mrm{Ext}}
\nc{\Fil}{\mrm{Fil}}
\nc{\Frob}{\mrm{Frob}}
\nc{\Gal}{\mrm{Gal}}
\nc{\GL}{\mrm{GL}}
\nc{\Hom}{\mrm{Hom}}
\nc{\hsr}{\mrm{H}}
\nc{\hpol}{\mrm{HP}}
\nc{\id}{\mrm{id}}
\nc{\im}{\mrm{im}}
\nc{\incl}{\mrm{incl}}
\nc{\length}{\mrm{length}}
\nc{\mchar}{\rm char\ }
\nc{\mpart}{\mrm{part}}
\nc{\ord}{\mrm{ord}}
\nc{\ql}{{\QQ_\ell}}
\nc{\qp}{{\QQ_p}}
\nc{\rank}{\mrm{rank}}
\nc{\rcot}{\mrm{cot}}
\nc{\rdef}{\mrm{def}}
\nc{\rdiv}{{\rm div}}
\nc{\rtf}{{\rm tf}}
\nc{\rtor}{{\rm tor}}
\nc{\res}{\mrm{res}}
\nc{\SL}{\mrm{SL}}
\nc{\Spec}{\mrm{Spec}}
\nc{\tor}{\mrm{tor}}
\nc{\Tr}{\mrm{Tr}}
\nc{\tr}{\mrm{tr}}

\nc{\Alg}{\mathbf{Alg}}
\nc{\Bax}{\mathbf{Bax}}
\nc{\bfk}{{\bf k}}
\nc{\bfone}{{\bf 1}}
\nc{\base}[1]{{a_{#1}}}
\nc{\detail}{\marginpar{\bf More detail}
    \noindent{\bf Need more detail!}
    \svp}
\nc{\Diff}{\mathbf{Diff}}   
\nc{\gap}{\marginpar{\bf Incomplete}\noindent{\bf Incomplete!!}
    \svp}
\nc{\FMod}{\mathbf{FMod}}
\nc{\Int}{\mathbf{Int}}
\nc{\Mon}{\mathbf{Mon}}
\nc{\remarks}{\noindent{\bf Remarks: }}
\nc{\Rep}{\mathbf{Rep}}
\nc{\Rings}{\mathbf{Rings}}
\nc{\Sets}{\mathbf{Sets}}

\nc{\dfootnote}[1]{{}}

\nc{\BA}{{\mathbb A}}
\nc{\CC}{{\mathbb C}}
\nc{\DD}{{\mathbb D}}
\nc{\EE}{{\mathbb E}}
\nc{\FF}{{\mathbb F}}
\nc{\GG}{{\mathbb G}}
\nc{\HH}{{\mathbb H}}
\nc{\LL}{{\mathbb L}}
\nc{\NN}{{\mathbb N}}
\nc{\QQ}{{\mathbb Q}}
\nc{\RR}{{\mathbb R}}
\nc{\TT}{{\mathbb T}}
\nc{\VV}{{\mathbb V}}
\nc{\ZZ}{{\mathbb Z}}

\nc{\cala}{{\cal A}}
\nc{\calc}{{\cal C}}
\nc{\cald}{\mathcal{D}}
\nc{\cale}{{\cal E}}
\nc{\calf}{{\cal F}}
\nc{\calg}{{\cal G}}
\nc{\calh}{{\cal H}}
\nc{\cali}{{\cal I}}
\nc{\call}{{\cal L}}
\nc{\calm}{{\cal M}}
\nc{\caln}{{\cal N}}
\nc{\calo}{{\cal O}}
\nc{\calp}{{\cal P}}
\nc{\calr}{{\cal R}}
\nc{\calt}{{\cal T}}
\nc{\calw}{{\cal W}}
\nc{\calx}{{\cal X}}
\nc{\CA}{\mathcal{A}}

\nc{\fraka}{{\mathfrak a}}
\nc{\frakA}{{\mathfrak A}}
\nc{\frakB}{{\mathfrak B}}
\nc{\frakm}{{\mathfrak m}}
\nc{\frakp}{{\mathfrak p}}
\nc{\frakS}{{\mathfrak S}}

\font\cyr=wncyr10 %wncyr10

\begin{document}
\title{Baxter Algebras and Differential Algebras
}
\author{Li Guo\\
Department of Mathematics and Computer Science, \\
Rutgers University at Newark, \\
Newark, NJ 07102, USA \\
E-mail: liguo@newark.rutgers.edu}
\date{}
\maketitle
% ----------------------------------------------------------------

\begin{center}
{\large In memory of Professor Chuan-Yan Hsiong
}
\end{center}

\medskip

\begin{abstract}
A Baxter algebra is a commutative algebra $A$ that carries a 
generalized integral operator. 
In the first part of this paper we review past work 
of Baxter, Miller, Rota and Cartier in this area and 
explain more recent work on explicit
constructions of free Baxter algebras that extended the constructions 
of Rota and Cartier. 
In the second part of the paper we will use these
explicit constructions to relate Baxter algebras to  
Hopf algebras and give applications of Baxter algebras 
to the umbral calculus
in combinatorics. 
\end{abstract}

\setcounter{section}{-1}

\section{Introduction}

This is a survey article on Baxter algebras, with emphasis on 
free Baxter algebras and their applications in probability theory, 
Hopf algebra and umbral calculus. This article can be read 
in conjunction with the excellent introductory article of 
Rota \cite{Ro2}. See also \cite{RS,Ro3}. 

\addsm

A Baxter algebra is a commutative algebra $R$ with a 
linear operator $P$ that satisfies the Baxter identity 
\begin{equation}
 P(x)P(y)=P(xP(y))+P(yP(x))+\lambda P(xy), \forall x, y\in R, 
\mlabel{eq:bax}
\end{equation}
where $\lambda$ is a pre-assigned constant, called the weight, 
from the base ring of $R$. 

\subsection{Relation with differential algebra}
The theory of Baxter algebras is related to differential algebra 
just as the integral analysis is related to 
the differential analysis. 

%\begin{enumerate}
%\item
Differential algebra originated from differential equations, 
while the study of Baxter 
algebras originated from the algebraic study by Baxter\cite{Ba} 
on integral equations arising from fluctuation theory in 
probability theory. 
%\item
Differential algebra has provided the motivation for some of the 
recent studies on Baxter algebras. 
A Baxter algebra of weight zero is 
an integration algebra, i.e., an algebra 
with an operator that satisfies the integration by parts formula 
(see Example 1 below).  
The motivation of the recent work in \cite{GK0,GK1,GK2} 
is to extend the beautiful theory of differential algebra~\cite{Ko} 
to integration algebras. 
%\item
On the other hand, the Baxter operators, regarded as a twisted family 
of integration operators, 
motivated the study of a twisted family of differential operators, 
generalizing the differential operator (when the weight is 0) and 
the difference operator (when the weight is 1). 
%\item
%Differential algebra and Baxter algebra have similar applications 
%to combinatorics and Hopf algebra. 
%\end{enumerate}

\subsection{Some history}
%\been
%\item
In the 1950's and early 1960's, several spectacular results were 
obtained in the fluctuation theory of sums of independent random 
variables by Anderson~\cite{An}, Baxter~\cite{Ba1}, 
Foata~\cite{Fo} and Spitzer~\cite{Sp}. 
The most important 
result is Spitzer's identity (see 
Proposition~\ref{pp:si}) which was applied to show
that certain functionals of sums of independent random 
variables, such as the maximum and the number of positive partial 
sum, were independent of the particular distribution. In an important paper \cite{Ba}, 
Baxter deduced Spitzer's identity and several other identities 
from identity (\ref{eq:bax}). This identity was further 
studied by Wendel \cite{We}, Kingman~\cite{Ki} and 
Atkinson~\cite{At}. 

%\item
Rota \cite{Ro1} realized the algebraic and combinatorial 
significance 
of this identity and started a systemic study of the algebraic 
structure of Baxter algebras. 
Free Baxter algebras were constructed by him \cite{Ro1,RS} and 
Cartier \cite{Ca}. 
%\item 
Baxter algebras were also applied to the study of Schur 
functions \cite{Th1,Mu,Wi}, hypergeometric functions and 
symmetric functions, and are closely related to several
areas in algebra and geometry, such as quantum groups and iterated
integrals, as well as differential algebra.
The two articles by Rota~\cite{Ro2,Ro3} include surveys in this 
area and further references. 
%See the articles~\cite{Ro2,Ro3} of Rota for surveys in 
%this area and for further references. 
%\item 
Rota's articles helped to revive the study of Baxter 
algebras in recent years: Baxter sequences in~\cite{Wi}, 
free Baxter algebras in~\cite{GK1,GK2,Gu1,Gu2} 
and applications in~\cite{AGKO,Gu3}
%\enen

Despite the close analogy between differential algebras and Baxter 
algebras, in particular integration algebras, 
relatively little is known about Baxter algebras in comparison 
with differential algebras. 
It is our hope that this article will further promote the study 
of Baxter algebras and related algebraic structures.

\subsection{Outline}

After introducing notations and examples in Section~\ref{s:dfn}, 
we will focus on the construction of free Baxter algebras 
in Section~\ref{s:free}. 
We will give three constructions of free Baxter algebras. 
We first explain Cartier's construction using brackets, followed 
by a similar construction using a generalization 
of shuffle products. 
These two constructions are ``external" in the
sense that each is a free Baxter algebra obtained without reference
to any other Baxter algebra. 
%These two constructions are exterior. 
We then explain Rota's 
standard Baxter algebra which chronologically came first. 
Rota's construction is an ``internal" construction,
obtained as a Baxter subalgebra inside a naturally
defined Baxter algebra whose construction traces back to 
Baxter~\cite{Ba}.

In Section~\ref{s:app}, we give two applications of Baxter algebras. 
We use free Baxter algebras to construct a new class of Hopf 
algebras, generalizing the classical divided power Hopf algebra. 
We then use Baxter algebras to give an interpretation and 
generalization of the umbral calculus. 
Other applications of free Baxter algebras can be found in 
Section~\ref{s:dfn}, relating Baxter operators to integration and 
summation, and in Section~\ref{s:free}, proving the famous formula 
of Spitzer. 
%\addsm

We are not able to include some other 
work on Baxter algebras, for example on Baxter sequences and 
the Young tableau~\cite{Th1,Wi}, 
and on zero divisors and chain conditions in 
free Baxter algebras \cite{Gu1,Gu2}. 
We refer the interested readers to the original literature. 
%\delmed

%\section{Baxter algebra}

\section{Definitions, examples and basic properties}
\mlabel{s:dfn}
%\delsm
\subsection{Definitions and examples}
We will only consider rings and algebras with identity 
in this paper. If $R$ is the ring or algebra, the identity 
will be denoted by $\bfone_R$, or by 1 if there is no 
danger of confusion. 
%\delsm
\begin{defn}
Let $C$ be a commutative ring. 
Fix a $\lambda$ in $C$. 
A {\em Baxter $C$-algebra {\rm (}of weight $\lambda${\rm )}}
is a commutative $C$-algebra $R$ together with a {\em Baxter 
operator {\rm (}of weight $\lambda${\rm )}} on $R$, that is, 
a $C$-linear operator $P:R\to R$
such that
$$ P(x)P(y)=P(xP(y))+P(yP(x))+\lambda P(xy),$$
for any $x,y\in R.$
\end{defn}
\addmed

Let $\Bax_C=\Bax_{C,\lambda}$ denote the category of Baxter $C$-algebras 
of weight $\lambda$ in which 
the morphisms are algebra homomorphisms that commute with 
the Baxter operators. 
%\delmed

%\subsection{Examples}
%\delsm
There are many examples of Baxter algebras. 
%\delsm
\begin{exam}({\bf Integration})
Let $R$ be $\Cont(\RR)$, the ring of continuous functions on $\RR$. 
For $f$ in $\Cont(\RR)$, define $P(f)\in \Cont(\RR)$ by 
\[ P(f)(x)=\int_0^x f(t)dt,\  x\in \RR. \]
Then $(\Cont(\RR),P)$ is a Baxter algebra of weight zero. 
\end{exam}

\begin{exam}({\bf Divided power algebra})
This is the algebra $$R=\bigoplus_{n\geq 0} C e_n$$ on which 
the multiplication is defined by 
$$e_m e_n =\bincc{m+n}{m} e_{m+n},\ m,\ n\geq 0.$$
The operator $P: R\to R$, where $P(e_n)=e_{n+1},\ n\geq 0,$ is 
a Baxter operator of weight zero. 
\label{ex2}
\end{exam}

\begin{exam}({\bf Hurwitz series)}
Let $R$ be $$ HC:=\{ (a_n) | a_n \in C, n\in \NN \},$$
the ring of Hurwitz series~\cite{Ke}. 
The addition is defined componentwise, and the multiplication is given
by $(a_n)(b_n) = (c_n)$, where 
$c_n = \sum_{k = 0}^n \bincc{n}{k}a_kb_{n-k}$. 
Define 
$$P : HC \rar HC,\ P((a_n)) =(a_{n-1}), {\rm\ where\ } 
    a_{-1}=0.$$ 
Then $HC$ is a Baxter algebra of weight $0$ which is  
the completion of the divided power algebra. 
\label{ex3}
\end{exam}
\addmed

We will return to Example~\ref{ex2} and \ref{ex3} 
in Section~\ref{ss:sha}. 

\begin{exam}({\bf Scalar multiplication}) Let $R$ be any $C$-algebra. 
For a given $\lambda\in C$, define 
$$P_\lambda:R \to R, 
x\mapsto -\lambda x, \forall\, x\in R.$$ 
Then $P_\lambda$ is a  Baxter operator of weight $\lambda$
on $R$.
\end{exam}

\begin{exam}({\bf Partial sums})
This is one of the first examples of a Baxter algebra, 
introduced by Baxter \cite{Ba}. Let $A$ be any $C$-algebra. 
Let 
$$R=\prod_{n\in \NN_+} A 
    =\{ (a_1,a_2,\ldots ) | a_n\in A,\ n\in \NN_+\}.$$
with addition, multiplication and scalar product defined 
entry by entry. Define $P:R\to R$ to be the ``partial sum" 
operator: 
$$P(a_1,a_2, \ldots) =\lambda 
(0, a_1, a_1+a_2,a_1+a_2+a_3, \ldots).$$
Then $P$ is a Baxter operator of weight $\lambda$ on $R$. 
\end{exam}
\addsm

We will return to this example in Section~\ref{ss:st}

\begin{exam}
({\bf Distributions~\cite{Ba}}) Let $R$ be the Banach 
algebra of functions 
$$ \varphi(t) =\int^\infty_{-\infty} e^{itx}dF(x) $$
where $F$ is a function such that 
$||\varphi||\defeq \int^\infty_{-\infty}|dF(x)|<\infty$ 
and such that $F(-\infty)\defeq \displaystyle{\lim_{x\to -\infty}} F(x)$ 
exists. The addition and multiplication are defined pointwise. 
Let $P(\varphi)(t)=\int^\infty_0 e^{itx}dF(x)+F(0)-F(-\infty)$. 
Then $(R,P)$ is a Baxter algebra of weight $-1$.
\label{ex6}
\end{exam}
\addsm

We will come back to this example in Proposition~\ref{pp:si}.

\begin{exam}
This is an important example in combinatorics~\cite{RS}. 
Let $R$ be the ring of functions 
$f: \RR\to \RR$ with finite support in which the product is 
the convolution: 
$$ (fg)(x):=\sum_{y\in \RR} f(y)g(x-y),\ x\in \RR.$$
Define $P: R\to R$ 
by 
$$ P(f)(x) =\sum_{y\in \RR, \max (0,y)=x} f(y),\ x\in \RR.$$
Then $(R,P)$ is a Baxter algebra of weight $-1$.
\end{exam}

\subsection{Integrations and summations}
Define a system of polynomials $\Phi_n(x)\in \QQ[x]$ 
by the generating function
$$ \frac{t(e^{xt}-1)}{e^t-1}=\sum_{n=0}^\infty 
\Phi_n(x)\frac{t^n}{n!}.$$
Then we have 
$$ \Phi_n(x)=B_n(x)-B_n,$$
where $B_n(x)$ is the $n$-th Bernoulli polynomial, defined 
by the generating function
$$\frac{te^{xt}}{e^t-1}=\sum_{n=0}^\infty 
B_n(x)\frac{t^n}{n!}\, ,$$
and $B_n=B_n(0)$ 
is the $n$-th Bernoulli number. It is well-known 
(see \cite[Section~15.1]{IR}) that 
\begin{equation}
 \Phi_{n+1}(k+1)=(n+1) \sum_{r=1}^k r^n
\mlabel{eq:ber}
\end{equation} 
for any integer $k\geq 1$. 
The following property of Baxter algebras is due to 
Miller \cite{Mi1}. 

\begin{prop}
\mlabel{pp:bp}
Let $C$ be a $\QQ$-algebra and let $R=C[t]$. 
Let $P$ be a $C$-linear operator on $R$ such that $P(1)=t$. 
\begin{enumerate}
\delsm\item \label{ber1}
The operator $P$ is a Baxter operator of weight $0$ if and only if 
$P(t^n)=\frac{1}{n+1} t^{n+1}$.
\delsm\item \label{ber2}
The operator $P$ is a Baxter operator of weight $-1$ if and only 
if $P(t^n)=\frac{1}{n+1}\Phi_{n+1}(t+1)$.
\delsm\item \label{ber3}
The operator $P$ is a Baxter operator of weight $-1$ if and only 
if $P(t^n)(k)=\displaystyle{\sum_{r=1}^k r^n}$ for every 
integer $k\geq 1$.
\end{enumerate}
\end{prop}
\delsm
\begin{proof}
\ref{ber1}: The only if part can be easily proved by induction. For details 
see the proof of Theorem 3 in \cite{Mi1}. The if part is obvious. 
\addsm

\ref{ber2} is Theorem 4 in \cite{Mi1}. 
\addsm

\ref{ber3}: Because of \ref{ber2}, we only need to show 
$$ P(t^n)=\frac{1}{n+1}\Phi_{n+1}(t+1) \Leftrightarrow 
    P(t^n)(k)=\sum_{r=1}^k r^n, \forall k\geq 1.$$
($\Rightarrow$) follows from (\ref{eq:ber}). 
($\Leftarrow$) can be seen easily, say from
\cite[Lemma 5]{Mi1}. 
\end{proof}

Using this proposition and free Baxter algebras that we will 
construct in the next section, we will prove 
the following property of Baxter algebras in Section~\ref{ss:sha}. 

\begin{prop}
\mlabel{pp:bp2}
Let $C$ be a $\QQ$-algebra and 
let $(R,P)$ be a Baxter $C$-algebra of weight $\lambda$. 
Let $t$ be $P(1)$.
\begin{enumerate}
\delsm\item
If $\lambda=0$, then 
$P(t^n)=\frac{1}{n+1} t^{n+1}$ for all $n\geq 1$.
\delsm\item
If $\lambda=-1$, then
$P(t^n)=\frac{1}{n+1}\Phi_{n+1}(t+1)$.
\end{enumerate}
\end{prop}
\addmed

Because of Proposition~\ref{pp:bp} and \ref{pp:bp2}, 
a Baxter operator of weight zero is 
also called an anti-derivation or integration, and a Baxter operator 
of weight $-1$  is also called a summation operator~\cite{Mi1}. 

\section{Free Baxter algebras}
\mlabel{s:free}
Free objects are usually defined to be generated by sets. 
We give the following more general definition. 
%\subsection{Definitions}
%\delsm
\begin{defn}
Let $A$ be a $C$-algebra. A Baxter $C$-algebra 
$(F_C(A),P_A)$, together with a $C$-algebra homomorphism $j_A:A\to F_C(A)$, 
is called a {\em free Baxter $C$-algebra on $A$ {\rm (}of weight 
$\lambda${\rm )}}, if, 
for any Baxter $C$-algebra $(R,P)$ of weight $\lambda$ 
and any $C$-algebra
homomorphism $\varphi:A\rar R$, there exists a unique Baxter
$C$-algebra homomorphism $\tilde{\varphi}:(F_C(A),P_A)\rar
(R,P)$ such that the  diagram
\[\xymatrix{
A \ar[rr]^(0.4){j_A} \ar[drr]_{\varphi}
    && F_C(A) \ar[d]^{\tilde{\varphi}} \\
&& R } \]
commutes.
\end{defn}
\addmed

Let $X$ be a set. One can define a free Baxter $C$-algebra $(F_C(X),P_X)$ 
on $X$ in a similar way. This Baxter algebra is naturally isomorphic to 
$(F_C(C[X]),P_{C[X]})$, where $C[X]$ is the polynomial algebra over $C$ 
generated by the set $X$. 

Free Baxter algebras play a central
role in the study of Baxter algebras. Even though the existence of
free Baxter algebras follows from the general theory of universal
algebras~\cite{Ma}, in order to get a good understanding of free Baxter
algebras, it is desirable to find concrete constructions. 
This is in analogy to the important role 
played by the ring of polynomials in the study of 
commutative algebra. 

Free Baxter algebras on sets were first constructed by 
Rota~\cite{Ro1} and Cartier~\cite{Ca} in the 
category of Baxter algebras with no identities 
(with some restrictions on the weight and the base ring).  
Two more general constructions have been obtained 
recently~\cite{GK1,GK2}.  
The first construction is in terms of {\it mixable shuffle 
products} which generalize the well-known 
{\it shuffle products} of path integrals 
developed by Chen~\cite{Ch} and Ree~\cite{Re}. 
The second construction is modified after the construction 
of Rota.

The shuffle product construction of a free Baxter algebra
has the advantage that its module
structure and Baxter operator can be easily described.
The construction of a free Baxter algebra as a standard Baxter 
algebra
has the advantage that its multiplication is very simple. 
There is a canonical isomorphism
between the shuffle Baxter algebra and
the standard Baxter algebra. 
This isomorphism enables us to make use of properties of 
both constructions.
%\newpage

\subsection{Free Baxter algebras of Cartier}
\mlabel{ss:ca}
We first recall the construction of free Baxter algebras 
by Cartier\cite{Ca,GK1}. The original construction of Cartier 
is for free objects on a set, in the category of algebras with 
no identity and with weight $-1$. We will extend his construction 
to free objects in the category of algebras with identity and with 
arbitrary weight. We will still consider free objects on sets. 
In order to consider free objects on a $C$-algebra $A$, it will 
be more convenient to use the tensor product notation to be 
introduced in Section~\ref{ss:sha}. 
%For further details, see\cite{GK1,Ca}. 

Let $X$ be a set. Let $\lambda$ be a fixed element in $C$. 
Let $M$ be the free commutative semigroup with identity on $X$.
Let $\widetilde{X}$ denote the set of symbols of the form
\[
 u_0\cdot [\ ],\ u_0\in M, \]
and
\[ u_0 \cdot [u_1,\ldots,u_m], m\geq 1,\ 
    u_0, u_1,\ldots,u_m\in M.
\]
Let $\frakB(X)$ be the free $C$-module on $\widetilde{X}$. 
Cartier gave a $C$-bilinear multiplication
$\shprc$ 
on $\frakB(X)$ by defining
\begin{eqnarray*}
(u_0 \cdot [\ ])\shprc (v_0\cdot [\ ])&=& u_0v_0 \cdot [\ ],\\
(u_0 \cdot [\ ])\shprc (v_0 \cdot [v_1,\ldots,v_n])
    &=& (v_0\cdot [v_1,\ldots,v_n])\shprc (u_0\cdot [\ ])\\
    &=& u_0 v_0\cdot [v_1,\ldots,v_n],
\end{eqnarray*}
and
\begin{eqnarray*}
\lefteqn{
(u_0 \cdot [u_1,\ldots,u_m])\shprc (v_0 \cdot [v_1,\ldots,v_n])}\\
    &=& \sum_{(k,P,Q)\in \overline{S}_c(m,n)}
    \lambda^{m+n-k} u_0 v_0 \cdot
    \Phi_{k,P,Q}([u_1,\ldots,u_m],[v_1\ldots,v_n]).
\end{eqnarray*}
Here $\overline{S}_c(m,n)$ is the set of triples $(k,P,Q)$ 
in which $k$ is an integer between $1$ and $m+n$, 
$P$ and $Q$ are ordered subsets of $\{1,\ldots, k\}$ with the
natural ordering such that
$P\cup Q=\{1,\ldots,k\}$, $\mid\!\! P\!\!\mid\ =m$ and
$\mid\!\! Q\!\!\mid\ =n$. 
For each $(k,P,Q)\in \overline{S}_c(m,n)$, 
$\Phi_{k,P,Q}([u_1,\ldots,u_m],[v_1,\ldots,v_n])$
is the element $[w_1,\ldots,w_k]$ in $\widetilde{X}$ defined by

\[
w_j=\left \{ \begin{array}{ll}
    u_\alpha, & {\rm\ if\ } j {\rm\ is\ the\ }
    \alpha\mbox{\rm -th\ element\ in\ } P {\rm\ and\ } j\not\in Q;\\
    v_\beta, & {\rm\ if\ } j {\rm\ is\ the\ }
    \beta\mbox{\rm -th\ element\ in\ } Q {\rm\ and\ } j\not\in P;\\
    u_\alpha v_\beta, & {\rm\ if\ } j {\rm\ is\ the\ }
    \alpha\mbox{\rm -th\ element\ in\ } P    
     {\rm\ and\ the\ } \beta\mbox{\rm -th\ element\ in\ } Q.
    \end{array} \right . \]
Define a $C$-linear operator $P_X^c$ on $\frakB(X)$ by 
\begin{eqnarray*}
P_X^c(u_0\cdot[\ ])&=&1\cdot [u_0], \\
P_X^c(u_0\cdot [u_1,\ldots,u_m])&=&1\cdot [u_0,u_1,\ldots,u_m].
\end{eqnarray*} 
The following theorem is a modification of Theorem 1 in 
\cite{Ca} and can be proved in the same way. 

\begin{thm} 
\mlabel{thm:ca}
The pair $(\frakB(X),P_X^c)$ is a
free Baxter algebra on $X$ of weight $\lambda$ 
in the category $\Bax_C$. 
\end{thm}
\delmed
\subsection{Mixable shuffle Baxter algebras}
\mlabel{ss:sha}
We now describe the mixable shuffle Baxter algebras. It gives 
another construction of free Baxter algebras. 
An advantage of this construction is that it is 
related to the well-known shuffle products. 
It also enables us to consider the free object on a 
$C$-algebra $A$. 

Intuitively, to form the shuffle product,  
one starts with two decks of cards and puts together all 
possible shuffles of the two decks. Similarly, 
to form the mixable shuffle product, one starts with 
two decks of charged cards, one deck positively charged
and the other negatively charged. 
When a shuffle of the two decks is taken, 
some of the adjacent 
pairs of cards with opposite charges are allowed to be merged into 
one card. When one puts all such ``mixable shuffles" 
together, with a proper measuring of the number of 
pairs that have been merged, one gets the mixable shuffle product. 
We now give a precise description of the construction. 

\subsubsection{Mixable shuffles}
For $m,n\in \NN_+$,
define the set of {\em $(m,n)$-shuffles} by
\addsm
\[ S(m,n)=
 \left \{ \sigma\in S_{m+n}
    \begin{array}{ll} {} \\ {} \end{array} \right .
\left |
\begin{array}{l}
\sigma^{-1}(1)<\sigma^{-1}(2)<\ldots<\sigma^{-1}(m),\\
\sigma^{-1}(m+1)<\sigma^{-1}(m+2)<\ldots<\sigma^{-1}(m+n)
\end{array}
\right \}.\]
Here $S_{m+n}$ is the symmetric group on $m+n$ letters. 
Given an $(m,n)$-shuffle $\sigma\in S(m,n)$,
a pair of indices $(k, k+1)$,\ $1\leq k< m+n$, is
called an {\em admissible pair} for $\sigma$
if $\sigma(k)\leq m<\sigma(k+1)$.
Denote $\calt^\sigma$ for the set of admissible pairs for $\sigma$.
For a subset $T$ of $\calt^\sigma$, call the pair
$(\sigma,T)$ a {\em mixable $(m,n)$-shuffle}.
Let $\mid\! T\!\mid$ be the cardinality of $T$.
Identify $(\sigma,T)$ with $\sigma$
if $T$ is the empty set.
Denote
\[ \bs (m,n)=\{ (\sigma,T)\mid \sigma\in S(m,n),\
    T\subset \calt^\sigma\} \]
for the set of {\em mixable $(m,n)$-shuffles}.
%Also denote
%\[ s(m,n)=\mid \bs(m,n)\mid .\]

%Intuitively, a $(m,n)$-shuffle is a permutation $\sigma$
%of $\{1,\ldots,m, m+1,\ldots, m+n\}$ such that
%the natural order of $\{1,\ldots,m\}$ and $\{m+1,\ldots,m+n\}$
%are preserved in
%$\{ \sigma(1),\ldots,\sigma(m),\sigma(m+1),\ldots,\sigma(m+n)\}$.
%Further, a mixable $(m,n)$-shuffle $(\sigma,T)$
%is a $(m,n)$-shuffle $\sigma$ in which pairs of indices
%from $T$ represent positions where ``merging" will occur.

\begin{exam}
There are three $(2,1)$ shuffles: 
$$\sigma_1=\left ( \begin{array}{ccc} 1 & 2& 3 \\
    1&2&3 \end{array} \right ) ,\,
\sigma_2=\left ( \begin{array}{ccc} 1 & 2& 3 \\
    1&3&2 \end{array} \right ),  \,
\sigma_3=\left ( \begin{array}{ccc} 1 & 2& 3 \\
    3&1&2 \end{array} \right ). $$
The pair $(2,3)$ is an admissible pair for $\sigma_1$.
The pair $(1,2)$ is an admissible pair for $\sigma_2$.
There are no admissible pairs for $\sigma_3$.
\end{exam}
\addmed

For $A\in \Alg_C$ and $n\geq 0$, let $A^{\otimes n}$ be the 
$n$-th tensor power of $A$ over $C$ with the convention 
$A^{\otimes 0}=C$. 
For  
$x=x_1\otimes\ldots\otimes x_m\in A^{\otimes m}$,
$y=y_1\otimes \ldots\otimes y_n\in A^{\otimes n}$
and $(\sigma,T)\in \overline{S}(m,n)$,
the element
\[  \sigma (x\otimes y) =u_{\sigma(1)}\otimes u_{\sigma(2)} \otimes
    \ldots \otimes u_{\sigma(m+n)}\in A^{\otimes (m+n)},\]
where
\[ u_k=\left \{ \begin{array}{ll}
    x_k,& \quad  1\leq k\leq m,\\
    y_{k-m}, & \quad m+1\leq k\leq m+n, \end{array}
    \right . \]
is called a {\em shuffle} of $x$ and $y$;
the element
\[ \sigma(x\otimes y; T)= u_{\sigma(1)}\hts u_{\sigma(2)} \hts
    \ldots \hts u_{\sigma(m+n)} \in A^{\otimes (m+n-\mid T\mid)},\]
where for each pair $(k,k+1)$, $1\leq k< m+n$,
\[ u_{\sigma(k)}\hts u_{\sigma(k+1)} =\left \{\begin{array}{ll}
    u_{\sigma(k)} u_{\sigma(k+1)},  &
    \quad    (k,k+1)\in T\\
    u_{\sigma(k)}\otimes u_{\sigma(k+1)}, &
    \quad (k,k+1) \not\in T,
    \end{array} \right . \]
is called a {\em mixable shuffle} of $x$ and $y$.

\begin{exam}
For $x=x_1\otimes x_2\in A^{\otimes 2}$ and $y=y_1\in A$, 
there are three shuffles of $x$ and $y$: 
\delsm
$$x_1\otimes x_2\otimes y_1,\, x_1\otimes y_1\otimes x_2,\, 
y_1\otimes x_1\otimes x_2.$$
For the mixable shuffles of $x$ and $y$, we have 
in addition, 
$$ x_1\otimes x_2 y_1,\, x_1 y_1 \otimes x_2.$$
\end{exam}
\addmed

Now fix $\lambda\in C$. 
Define, for $x$ and $y$ above, 
\begin{equation}
x \shpr\!\!^+ y\
=\sum_{(\sigma,T)\in \bs (m,n)} \lambda^{\mid T\mid } \sigma(x\otimes y;T) 
    \in \bigoplus_{k\leq m+n} A^{\otimes k}.
\mlabel{eq:shuf}
\end{equation}
\delmed
\delmed
\begin{exam}
{\rm For $x,\ y$ in our previous example, 
$$x\shpr\!\!^+ y 
= x_1\otimes x_2\otimes y_1+ x_1\otimes y_1\otimes x_2
+y_1\otimes x_1\otimes x_2+
\lambda( x_1\otimes x_2 y_1+ x_1 y_1 \otimes x_2).$$
}
\end{exam}
\delsm

\subsubsection{Shuffle Baxter algebras}
The operation $\shpr\!^+ $ extends to a map
\[ \shpr\!^+\ : A^{\otimes m}\times A^{\otimes n} \rar
    \bigoplus_{k\leq m+n} A^{\otimes k},\, m,\, n\in \NN 
\delmed \delmed\]
by $C$-linearity.
Let
\[ \sha_C^+(A)=\sha_C^+(A,\lambda)= \bigoplus_{k\in\NN} A^{\otimes k}
=C\oplus A\oplus A^{\otimes 2}\oplus \ldots. \]
Extending by additivity, the binary operation
$\shpr^+$ gives a $C$-bilinear map 
\[ \shpr \!^+\ : \sha_C^+(A) \times \sha_C^+(A) \rar \sha_C^+(A) \]
with the convention that
\[ C\times A^{\otimes m} \rar A^{\otimes m} \]
is the scalar multiplication. 

\begin{thm} {\rm \cite{GK1}}
The mixable
shuffle product $\shpr \!^+$ defines 
an associative, commutative binary operation on
$\sha_C^+(A)= \bigoplus_{k\in\NN} A^{\otimes k}$, making it into
a $C$-algebra
with the identity $\bfone_C\in C=A^{\otimes 0}$.
\end{thm}

\addmed

Define 
$ \sha_C(A)=\sha_C(A,\lambda)=A\otimes_C \sha_C^+(A)$
to be the tensor product algebra. 
Define a $C$-linear endomorphism $P_A$ on
$\sha_C(A)$ by assigning
\[ P_A( x_0\otimes x_1\otimes \ldots \otimes x_n)
=\bfone_A\otimes x_0\otimes x_1\otimes \ldots\otimes x_n, \]
for all
$x_0\otimes x_1\otimes \ldots\otimes x_n\in A^{\otimes (n+1)}$
and extending by additivity.
Let $j_A:A\rar \sha_C(A)$ be the canonical inclusion map. 
Call $(\sha_C(A),P_A)$ the {\em {\rm (}mixable{\rm )} shuffle Baxter $C$-algebra on
$A$ of weight $\lambda$}. 
\addsm

For a given set $X$, we also let $(\sha_C(X),P_X)$ denote the
shuffle Baxter $C$-algebra $(\sha_C(C[X]),P_{C[X]})$, called the
{\em {\rm (}mixable{\rm )} shuffle Baxter $C$-algebra on $X$ 
{\rm(}of weight $\lambda${\rm )}.} Let
$j_X:X\to \sha_C(X)$ be the canonical inclusion map.

\begin{thm}
\mlabel{thm:shua} {\rm \cite{GK1}}
The shuffle Baxter algebra $(\sha_C(A),P_A)$, 
together with the natural embedding $j_A$, is a
free Baxter $C$-algebra on $A$ of weight $\lambda$. 
Similarly,
$(\sha_C(X),P_X)$, together with the natural embedding
$j_X$,  is a free Baxter $C$-algebra on $X$
of weight $\lambda$.
\end{thm}
\subsubsection{Relation with Cartier's construction}
The mixable shuffle product construction of free Baxter algebras 
is canonically isomorphic to Cartier's construction. 
Using the notations from Section~\ref{ss:ca}, 
we define a map $f: \widetilde{X}\to \sha_C(X)$ by
\begin{eqnarray*}
f(u_0\cdot [\ ])& =& u_0;\\
f(u_0\cdot [u_1,\ldots,u_m])
&=& u_0\otimes u_1\otimes\ldots\otimes u_m, 
\end{eqnarray*}
and extend it by $C$-linearity to a $C$-linear map
\[ f: \frakB(X)\to \sha_C(X). \]
The same argument as for Proposition 5.1 in \cite[Prop. 5.1]{GK1} 
can be used to prove the following
\begin{prop}
\mlabel{prop:inj}
$f$ is an isomorphism in $\Bax_C$.
\end{prop}
\addmed
Let $A$ be a $C$-algebra. 
Using the tensor product notation in the construction of 
$\sha_C(A)$, one can extend Cartier's construction 
(Theorem~\ref{thm:ca}) and Proposition~\ref{prop:inj} for 
free Baxter algebras on $A$. 
\addmed
%\newpage

\subsubsection{Special cases}
We next consider some special cases of the shuffle Baxter 
algebras. 
\addmed

\noindent
{\bf Case 1: $\lambda=0$.} 
In this case, 
$\sha_C(A)$ is the usual shuffle algebra generated 
by the $C$-module $A$. It played a central role 
in the work of K.T. Chen~\cite{Ch} on path integrals and  
is related to many areas of pure and applied mathematics. 
\addsm

\noindent
{\bf Case 2: $X=\phi$\,.}
Taking $A=C$, we get
\[ \sha_C(C)=\bigoplus_{n=0}^\infty C^{\otimes (n+1)}
= \bigoplus_{n=0}^\infty C \bfone^{\otimes (n+1)}, \]
where
$\bfone^{\otimes (n+1)}
= \underbrace{\bfone_C \otimes \ldots \otimes \bfone_C}
_{(n+1)-{\rm factors}}$. In this case the mixable shuffle
product formula~(\ref{eq:shuf}) gives

\begin{prop}
%\margin{prop:unit}
\label{prop:unit}
For any $m,n\in \NN$,
\[ \bfone^{\otimes (m+1)}  \bfone^{\otimes (n+1)} =
\sum_{k=0}^m \binc{m+n-k}{n}\binc{n}{k} \lambda^k
\bfone^{\otimes (m+n+1-k)}.\]
\end{prop}
\addmed

We are now ready to prove Proposition~\ref{pp:bp2}\,: 
By Proposition 4.2 in \cite{Gu3}, the map 
\begin{eqnarray*}
 \sha_C(C) & \to &C[x],\\ 
   \bfone^{\otimes (n+1)} & \mapsto &
\frac{x(x-\lambda) \cdots (x-\lambda(n-1))}{n!}   
    ,\ n\geq 0, \addsm
\end{eqnarray*}
is an isomorphism of $C$-algebras. 
Then $P_C$ enables us to define a Baxter operator $Q$ on 
$C[x]$ through this isomorphism and  
we have $Q(1)=x$. Then by Proposition~\ref{pp:bp}, we have 
$$ Q(x^n) = \left \{ \begin{array}{ll} 
    \frac{1}{n+1} x^{n+1},& {\rm\ if\ } \lambda=0,\addmed\\
    \frac{1}{n+1} \Phi_{n+1}(x+1), &{\rm\ if\ } \lambda=-1.
\end{array} \right . $$
Now let $(R,P)$ be any Baxter $C$-algebra. 
By the universal property of $(C[x],Q)$ ($\cong (\sha_C(C),P_C)$) 
stated in Theorem~\ref{thm:shua}, 
there is a unique homomorphism 
$\widetilde{\varphi}: (C[x],Q) \to (R,P)$ of Baxter algebras such that 
$\widetilde{\varphi}(x)=P(1)$. 
Let $t = P(1)=\widetilde{\varphi}(Q(1))=\widetilde{\varphi}(x)$ 
and since $\Phi_{n+1}(x+1)$ has 
coefficients in $\QQ\subset C$, we have 
$$ P(t^n) = \widetilde{\varphi}(Q(x^n))= 
\widetilde{\varphi}\left(\frac{1}{n+1} x^{n+1}\right) 
=\frac{1}{n+1} t^{n+1}$$ when 
$\lambda=0$ and 
$$ P(t^n) =\widetilde{\varphi}(Q(x^n))= 
\widetilde{\varphi}\left(\frac{1}{n+1} \Phi_{n+1}(x+1)\right)
=\frac{1}{n+1} \Phi_{n+1}(t+1)$$
when $\lambda=-1$.
This proves Proposition~\ref{pp:bp2}. 
%\proofend
\addmed

\noindent
{\bf Case 3: $\lambda=0$ and $X=\phi$.}
%: divided power algebra}
Taking the ``pull-back" of Cases 1 and 2, we get 
\begin{equation}
\begin{CD}
    \{ \sha_C(\phi, 0)\} @>{\subset}>> 
     \{\sha_C(X,0) | X\in \Sets \} \\
@V\cap VV  @ VV \cap V\\
\{ \sha_C(\phi,\lambda) | \lambda \in C \} @>{\subset}>>
 \{\sha_C(X,\lambda) | X\in \Sets, \lambda \in C\}
\end{CD}
\addmed
\mlabel{eq:free}
\end{equation}
Thus $\sha_C(\phi, 0)$ is the {\em divided power algebra} 
\[ \sha_C(\phi, 0)= \bigoplus_{k\in \NN} C e_k, \ 
e_n e_m =\binc{m+n}{m} e_{m+n}\]
in Example~\ref{ex2}. 

%\newpage

\subsubsection{Variation: Complete shuffle Baxter algebras}
We now consider the completion of $\sha_C(A)$.

\addsm
Given $k\in \NN$,
$\Fil^k \sha_C(A)\defeq \bigoplus_{n\geq k} A^{\otimes (n+1)}$
is a Baxter ideal of $\sha_C(A)$.
Consider the infinite product of $C$-modules
$\widehat{\sha}_C(A)=\prod_{k\in \NN} A^{\otimes (k+1)}$.
It contains $\sha_C(A)$ as a dense subset with respect to the
topology defined by the filtration $\{\Fil^k \sha_C(A)\}$.
All operations of the Baxter $C$-algebra $\sha_C(A)$ are continuous
with respect to this topology. Hence they extend uniquely to
operations on $\widehat{\sha}_C(A)$,
making $\widehat{\sha}_C(A)$ a Baxter algebra
of weight $\lambda$, with the Baxter operator denoted by $\widehat{P}$.
It is called the {\em complete shuffle Baxter algebra on $A$}.
It naturally contains $\sha_C(A)$ as a Baxter subalgebra
and is a free object in the category of Baxter algebras
that are complete with respect to a canonical filtration
defined by the Baxter operator~\cite{GK2}.

When $A=C$ and $\lambda=0$, we have 
\[\widehat{\sha}_C(C,\lambda)= \prod_{k\in \NN} C e_k 
    \cong HC, \] 
the ring of Hurwitz series in Example~\ref{ex3}.

\subsection{Standard Baxter algebras}
\mlabel{ss:st}
The standard Baxter algebra constructed by Rota in~\cite{Ro1}
is a free object in
the category $\Bax_C^0$ of Baxter algebras not necessarily
having an identity. It is described as a Baxter subalgebra
of another Baxter algebra whose construction goes back to
Baxter~\cite{Ba}.
In Rota's construction, there are further restrictions 
that $C$ be a field of characteristic zero,
the free Baxter algebra obtained be
on a finite set $X$, and the weight $\lambda$ be $1$.
%Among these restrictions, $C$ being a field of characteristic zero
%and the weight being $1$ are easy to remove,
%as long as the weight is not a zero divisor, 
%and the same argument given by Rota can be carried through.
%The relaxation of other restrictions is a less trivial matter. 
%We will do this by making use of the shuffle Baxter algebra. 
By making use of shuffle Baxter algebras, 
we will show that Rota's description can be modified to yield 
a free Baxter algebra on an algebra in the category $\Bax_C$
of Baxter algebras with an identity,
with a mild restriction on the weight $\lambda$. 
We can also provide a similar construction for
algebras not necessarily having an identity,
and for complete Baxter algebras, but we will not explain it 
here. See~\cite{GK2}. 

We will first present Rota's construction, modified to give 
free objects in the category of Baxter algebras with identity. 
We then give the general construction. 

%This class of free Baxter algebras will be called Baxter
%algebras of Rota, or the internal free Baxter algebras,
%while the shuffle Baxter algebras will sometimes be called the external
%free Baxter algebras. 
%We will give a natural identification of
%the shuffle Baxter algebra and the Baxter algebra of Rota. 

\subsubsection{The standard Baxter algebra of Rota}
For details, see~\cite{Ro1,RS}.
%The original construction of was given for Baxter algebras
%over a field of characteristic zero.
%We will show that the
%same construction works for Baxter algebras over any commutative
%ring, provided the weight is not a zero divisor in
%the base ring.

As before, let $C$ be a commutative ring with an identity, 
and fix a $\lambda$ in $C$. Let $X$ be a given set.
For each $x\in X$, let $t^{(x)}$ be a sequence
$( t^{(x)}_1,\ldots, t^{(x)}_n,\ldots )$
of distinct symbols $t^{(x)}_n$.
We also require that the sets $\{t^{(x_1)}_n\}_n$ and
$\{t^{(x_2)}_n\}_n$ be disjoint for $x_1\neq x_2$ in $X$.
Denote
\[ \overline{X} = \bigcup_{x\in X} \{t^{(x)}_n \mid n\in \NN_+\} \]
and denote by $\frakA (X)$ the ring of sequences with entries in 
$C[\overline{X}]$, the $C$-algebra of polynomials with variables in
$\overline{X}$.
Thus the addition, multiplication and scalar multiplication by
$C[\overline{X}]$ in $\frakA(X)$ are defined componentwise.
Alternatively, for $k\in \NN_+$, denote $\gamma_k$ for the sequence
$(\delta_{n,k})_n$, where $\delta_{n,k}$ is the Kronecker delta. 
Then we can identify a sequence $(a_n)_n$ in $\frakA(X)$
with a series
\[ \sum_{n=1}^\infty a_n \gamma_n = a_1 \gamma_1 +a_2 \gamma_2 + \ldots .\]
Then the addition, multiplication and scalar multiplication by
$C[\overline{X}]$ are given termwise. 

Define
\[P_X^r=P_{X,\lambda}^r: \frakA(X)\to \frakA(X)\]
by
\[ P_X^r(a_1,a_2,a_3,\ldots)
=\lambda (0,a_1,a_1+a_2,a_1+a_2+a_3,\ldots).\]
In other words, each entry of $P_X^r (a),\ a=(a_1,a_2,\ldots),$
is $\lambda$ times the sum of the previous entries of $a$.
If elements in $\frakA(X)$ are described by series
$\sum_{n=1}^\infty a_n \gamma_n$ given above, then we simply have
\[ P_X^r\left(\sum_{n=1}^\infty a_n \gamma_n\right)=
    \lambda \sum_{n=1}^\infty \left(\sum_{i=1}^{n-1} a_i\right) \gamma_n.\]
It is well-known~\cite{Ba,Ro1} that, for $\lambda=1$,
$P_X^r$ defines a Baxter operator of weight $1$ on $\frakA(X)$.
It follows that, for any $\lambda\in C$,
$P_X^r$ defines a Baxter operator of weight $\lambda$ on $\frakA(X)$. 
Hence $(\frakA(X),P^r_X)$ is in $\Bax_C$. 

\begin{defn}
The {\em standard Baxter algebra} on $X$ is the Baxter subalgebra 
$\frakS(X)$ of $\frakA(X)$ generated
by the sequences $t^{(x)}=(t^{(x)}_1,\ldots,x^{(x)}_n,\ldots),\ x\in X$.
\end{defn}
\addmed

An important result of Rota~\cite{Ro1,RS} is 
\begin{thm}
\mlabel{thm:rota}
$(\frakS(X),P_X^r)$ is a free Baxter algebra on $X$ in the
category $\Bax_C$.
\end{thm}

\subsubsection{The standard Baxter algebra in general}
Given $A\in \Alg_C$, 
we now give an alternative construction of
a free Baxter algebra on $A$ in the category $\Bax_C$.

For each $n\in \NN_+$, denote by $A^{\otimes n}$ the $n$-th 
tensor power
algebra where the tensor product is taken over $C$.
Note that the multiplication on $A^{\otimes n}$ here is
different from the multiplication on $A^{\otimes n}$ when it
is regarded as a $C$-submodule of $\sha_C(A)$.

Consider the direct limit algebra 
\[ \overline{A}= \dirlim A^{\otimes n} \]
where the transition map is given by
\[ A^{\otimes n} \la A^{\otimes (n+1)},\
    x \mapsto x\otimes \bfone_A.\]
%To distinguish between the two contexts, we will use the notation
%$\sha_C^n(A)$ for $A^{\otimes (n+1)}\subseteq \sha_C(A)$. 
Let $\frakA(A)$ be the set of sequences with entries in
$\overline{A}$. Thus we have
\[ \frakA(A) = \prod_{n=1}^\infty \overline{A} \gamma_n
    =\left \{ \sum_{n=1}^\infty a_n \gamma_n,\ 
    a_n \in \overline{A} \right \}.\]
Define addition, multiplication and scalar multiplication
on $\frakA(A)$ componentwise, making $\frakA(A)$ into a
$\overline{A}$-algebra,
with the sequence $(1,1,\ldots)$ as the identity. 
Define
\[P_A^r=P_{A,\lambda}^r: \frakA(A)\to \frakA(A)\]
by
\[ P_A^r(a_1,a_2,a_3,\ldots)
=\lambda (0,a_1,a_1+a_2,a_1+a_2+a_3,\ldots).\]
Then $(\frakA(A),P^r_A)$ is in $\Bax_C$. 
For each $a\in A$, define $t^{(a)}=(t^{(a)}_k)_k$ in $\frakA(A)$ by 
\[ t^{(a)}_k=\otimes_{i=1}^k a_i
    (=\otimes_{i=1}^\infty a_i),\ a_i = \left \{
    \begin{array}{ll} a, & i=k\, ,\smallskip\\
    1, & i\neq k\,. \end{array} \right.  \]
\delmed\delmed
\begin{defn}
The {\em standard Baxter algebra} on $A$ is the Baxter subalgebra 
$\frakS(A)$ of $\frakA(A)$ generated
by the sequences
$t^{(a)}=(t^{(a)}_1,\ldots,t^{(a)}_n,\ldots),\ a\in A$.
%Let $P_X^R$ be the restriction of $P_X^r$ on $\frakS(X)$.
\end{defn}
\addmed

Since $\sha_C(A)$ is a free Baxter algebra on $A$,
the $C$-algebra morphism
\[ A \to \frakA(A), a\mapsto t^{(a)}\]
extends uniquely to a morphism in $\Bax_C$
\[ \Phi: \sha_C(A) \to \frakA(A). \]
\delbig\delmed

\begin{thm}
\mlabel{thm:s-r} {\rm \cite{GK2}}
Assume that $\lambda\in C$ is not a zero divisor in $\overline{A}$. 
The morphism in $\Bax_C$ 
\delmed\[ \Phi: \sha_C(A)\to \frakS(A) \]
induced by sending
$a\in A$ to $t^{(a)}=(t^{(a)}_1,\ldots,t^{(a)}_n,\ldots)$
is an isomorphism.
\end{thm}
\addmed

Consequently, when $\lambda$ is not a zero divisor in 
$\overline{A}$, 
$(\frakS(A),P_A^r)$ is a free Baxter algebra on $A$ in the
category $\Bax_C$. 

\subsubsection{Spitzer's identity}
As an application of the standard Baxter algebra, we recall 
the proof of Spitzer's identity by Rota~\cite{Ro1,RS}. 
Spitzer's identity is regarded as a remarkable stepping stone in 
the theory of sums of independent random variables and motivates  
of Baxter's identity. 
For other proofs of Spitzer's identity, see \cite{Ba,Ki,We,At,Ca}. 
We first present an algebraic formulation. 
\begin{prop} {\rm \cite{RS}}
\mlabel{pp:si}
Let $C$ be a $\QQ$-algebra. 
Let $(R,P)$ be a Baxter $C$-algebra of weight 1. Then for $b\in R$, 
we have 
\begin{equation}
 \exp\left (-P(\log(1+tb)^{-1}) \right ) 
    =\sum_{n=0}^\infty t^n (Pb)^{[n]}
\mlabel{eq:si1}
\end{equation}
in the ring of power series $R[[t]]$. 
Here $$
(Pb)^{[n]}=\underbrace{P(b(P(b \ldots (Pb) \ldots )))}_{n\mbox{\rm 
-}{\rm iteration}}$$
with the convention that $(Pb)^{[1]}=P(b)$ and $(Pb)^{[0]}=1$.
\end{prop}
\begin{proof}
First let $x=(x_1,x_2,\ldots)$ where $x_i,\ i\geq 1$, are symbols 
and let $X=\{x\}$. 
Consider the standard Baxter algebra
$\frakS(X)$. It is easy to verify that 
$$(Px)^{[n]}=(0,e_n(x_1),e_n(x_1,x_2),e_n(x_1,x_2,x_3),\ldots)$$
where $e_n(x_1,\ldots,x_m)$ is the elementary symmetric function 
of degree $n$ in the variables $x_1,\ldots,x_m$. By definition, 
$$P(x^k)=(0,x_1^k,x_1^k+x_2^k,x_1^k+x_2^k+x_3^k,\ldots,
p_k(x_1,\ldots,x_m),\ldots),$$
where $p_k(x_1,\ldots,x_m)=x_1^k+x_2^k+\ldots,x_m^k$ is the power 
sum symmetric function of degree $k$ in the variables 
$x_1,\ldots,x_m$. These two classes of symmetric functions are 
related by the well-known Waring's formula~\cite{RS}
$$ \exp\left (-\sum_{k=1}^\infty (-1)^kt^k p_k(x_1,\ldots,x_m)/k 
    \right) = \sum_{n=0}^\infty e_n(x_1,\ldots,x_m)t^n,\ 
\forall\ m\geq 1.$$
This proves 
\begin{equation}
 \exp\left (-P_X^r(\log(1+tx)^{-1})\right )
    =\sum_{n=0}^\infty t^n (P_X^r x)^{[n]}.
\mlabel{eq:si2}
\end{equation}

Next let $(R,P)$ be any Baxter $C$-algebra and let $b$ be any 
element in $R$. By the universal property of the free Baxter 
algebra $(\frakS,P_X^r)$, there is a unique Baxter algebra 
homomorphism $\tilde{\varphi}: \frakS\to R$ such that 
$\tilde{\varphi}(x)=b$. Since all the coefficients in the 
expansion of $\log (1+u)$ and $\exp (u) $ are rational, 
and $\tilde{\varphi}\circ P_X^r=P\circ \tilde{\varphi}$, 
applying $\tilde{\varphi}$ to 
(\ref{eq:si2}) gives the desired equation. 
\end{proof}

We can now specialize to the original identity of Spitzer, 
following Baxter~\cite{Ba} and Rota~\cite{RS}. 
Consider the Baxter algebra $(R,P)$ in Example~\ref{ex6}. 
Let $\{X_k\}$ be a sequence of independent random 
variables with identical distribution function $F(x)$ and 
characteristic function 
\addsm$$\psi(s)=\int_{-\infty}^\infty e^{isx} dF(x). 
\addsm $$
Let $S_n=X_1+\ldots+X_n$ and let 
$M_n=\max (0,S_1,S_2,\ldots,S_n)$. 
Let $F_n(x)={\rm Prob}(M_n<x)$ (${\rm Prob}$ for probability) be the distribution 
function of $M_n$. 
We note that, if $f(s)$ is the characteristic function of the 
random variable of $X$, 
then $P(f)(s)$ is the characteristic function of the random 
variable $\max(0,X)$. 
Applying Proposition~\ref{pp:si} to $b=\psi(s)$, we obtain
the identity first obtained by Spitzer:
\[
 \sum_{n=0}^\infty \int_0^\infty e^{isx} d F_n(x) 
    =\exp\left ( \sum_{k=1}^\infty \left ( \int_0^\infty e^{isx} 
    dF(x)+F(0)\right ) \right ).
%\mlabel{eq:si3}
\]

We refer the reader to \cite{RS} for the application of the 
standard Baxter algebra to the proof of some other identities, 
such as the Bohnenblust-Spitzer formula.

\section{Further applications of free Baxter algebras}
\mlabel{s:app}

\subsection{Overview}
Recall that the free Baxter algebra $\sha_C(A,\lambda)$ 
in the special case when $A=C$ and $\lambda=0$ 
is the divided power algebra. 
The divided power algebra and its completion 
are known to be related to  
\delmed
\begin{itemize}
\item 
crystalline cohomology and rings of $p$-adic periods 
in number theory, 
%\footnote{well-known}
\delmed
\item
shuffle products in differential geometry and topology, 
%\footnote{classic}
\delmed\item
Hopf algebra in commutative algebra,
%\footnote{well-known}
\delmed\item
Hurwitz series in differential algebra, 
%\mnote{cofree differential algebra}
\delmed\item
umbral calculus in combinatorics, and 
%\footnote{will explain later}
\delmed\item
incidence algebra in graph theory.
\delsm
\end{itemize}
By the ``pull-back" diagram (\ref{eq:free}), 
free (complete) Baxter algebras 
give a vast generalization of the (complete) divided power 
algebra, and so suggest a framework in which these connections 
and applications of the divided power algebra can   
be extended. 
We give two such connections and applications in 
the next two sections, one to Hopf algebras  
(Section~ \ref{ss:AGKO}) and one 
to the umbral calculus in combinatorics (Section~ \ref{ss:um}). 
\delmed

\subsection{Hopf algebra}
\mlabel{ss:AGKO}

\subsubsection{Definition of Hopf algebra} 
We recall some basic definitions and facts.
Recall that a {\it cocommutative 
$C$-coalgebra} is a triple $(A,\Delta,\varepsilon)$ 
where $A$ is a $C$-module, and $\Delta: A\to A\otimes A$ 
and $\varepsilon: A\to C$ are $C$-linear maps that make  
the following diagrams commute. 
\begin{equation}
\begin{array}{ccc} 
    A & \ola{\Delta} & A\otimes A \\
    \dap{\Delta} && \dap{\id\otimes \Delta}\\
    A\otimes A &\ola{\Delta\otimes \id} & A\otimes A\otimes A
\end{array}
\mlabel{coass}
\end{equation}
\begin{equation}
\begin{array}{ccccc}
C\otimes A& \stackrel{\varepsilon\otimes \id}{\longleftarrow} & 
A\otimes A & \ola{\id\otimes \varepsilon} & A\otimes C \\
&  {}_{\cong}\!\!\nwarrow &\uap{\Delta} & \nearrow_{\cong} \\
&& A&&
\end{array}
\mlabel{coun}
\end{equation}
\begin{equation}
\begin{array}{ccccc}
& & A & & \\
& {}^{\Delta} \!\!\swarrow & & \searrow^{\Delta}& \\
A\otimes A & & \ola{\tau_{A,A}} & &A\otimes A
\end{array}
\mlabel{cocom}
\end{equation}
where $\tau_{A,A}: A\otimes A\to A\otimes A$ is defined 
by $\tau_{A,A}(x\otimes y)=y\otimes x$. 

Recall that a $C$-{\it bialgebra} is a quintuple 
$(A,\mu,\eta,\Delta,\varepsilon)$ where $(A,\mu,\eta)$ 
is a $C$-algebra and $(A,\Delta,\varepsilon)$ is a $C$-coalgebra 
such that $\mu$ and $\eta$ are morphisms of coalgebras. 

Let $(A,\mu,\eta,\Delta,\varepsilon)$ be a $C$-bialgebra. 
For $C$-linear maps $f,\, g:A\to A$, the convolution 
$f\star g$ of $f$ and $g$ is the composition of the maps 
\[ A\ola{\Delta} A\otimes A \ola{f\otimes g} A\otimes A 
    \ola{\mu} A.\]
A $C$-linear endomorphism $S$ of $A$ is called an {\it antipode} 
for $A$ if 
\begin{equation}
 S\star \id_A = \id_A\star S =\eta\circ \varepsilon. 
\mlabel{anti}
\end{equation}
A {\it  Hopf algebra} is a bialgebra $A$ with an antipode $S$. 

%\newpage

\subsubsection{The main theorem}
On the Baxter algebra 
$\sha_C(C,\lambda)$, let $\mu$ be the canonical multiplication 
and let 
$\eta: C\hookrightarrow \sha_C(C,\lambda)$ be the unit map. 
Define a comultiplication $\Delta$, a counit $\varepsilon$ and 
an antipode $S$ by
\begin{eqnarray*}
& \ \ \  \Delta= \Delta_\lambda:& 
\sha_C(C,\lambda)\to \sha_C(C,\lambda) \otimes 
    \sha_C(C,\lambda),  \\
&&  \ \     \base{n}\mapsto  
\sum_{k=0}^{n} \sum_{i=0}^{n-k} (-\lambda)^k \base{i}\otimes 
    \base{n-k-i},\\
& \ \ \ \ \varepsilon=\varepsilon_\lambda: & 
\sha_C(C,\lambda)\to C,\ \base{n} \mapsto \left \{ 
\begin{array}{ll}
\bfone, & n=0,\\
\lambda\bfone, & n=1,\\
0, & n\geq 2 ,
\end{array} \right . \\
&\ \ \ \  S=S_\lambda: &\sha_C(C,\lambda)
\to \sha_C(C,\lambda),\ 
 \ \ \base{n} \mapsto 
  (-1)^n \sum_{v=0}^n \bincc{n-3}{v-3} \lambda^{n-v} \base{v}.
\label{S}
\end{eqnarray*}
The following result is proved in {\rm \cite{AGKO}}.
\begin{thm}The sextuple $(\sha_C(C,\lambda), \mu,\eta,\Delta,\varepsilon,S)$ 
is a Hopf $C$-algebra. 
\end{thm}

\subsection{The umbral calculus}
\mlabel{ss:um}
\subsubsection{Definition and examples}
For simplicity, we assume that $C$ is a $\QQ$-algebra
for the rest of the paper. 

The {\em umbral calculus} is the study and application of 
{\em polynomial sequences of binomial type}, i.e.,  
polynomial sequences $\{ p_n(x)\,|\, n\in \NN\}$ 
in $C[x]$ such that  
$$ p_n(x+y)=\sum_{k=0}^n \bincc{n}{k} p_k(x) p_{n-k}(y)
$$
in $C[x,y]$ for all $n$. 
Such a sequence behaves as if its terms  
are powers of $x$ and has 
found applications in several areas of pure and applied 
mathematics, including number theory and combinatorics, 
since the 19th century. 
There are many well-known sequences of binomial types.  

\begin{exams}
{\rm 
\been
\item
Monomials.  $x^n$. 
\item
 Lower factorial polynomials. 
$(x)_n=x(x-1)\cdots (x-n+1)$.
\item
 Exponential polynomials. 
$\phi_n(x)=\sum_{k=0}^n S(n,k) x^k$, where 
$S(n,k)$ with $ n, k\geq 0$ are the Stirling numbers of the 
second kind. 
\item
 Abel polynomials. Fix $a\neq 0$.  
$A_n(x)=x(x-an)^{n-1}$.
\item Mittag-Leffler polynomials. 
$M_n(x)=\sum_{k=0}^n \bincc{n}{k}\bincc{n-1}{n-k} 
    2^k (x)_k$.
\item
Bessel polynomials. 
$y_n(x)=\sum_{k=0}^n \frac{(n+k)!}{(n-k)!k!} 
    (\frac{x}{2})^k$ 
(a solution to the Bessel 
equation
$x^2 y'' +(2x+2)y'+n(n+1) y=0).$
\enen
}
\end{exams}
\addsm
There are also Bell polynomials, Hermite polynomials, Bernoulli 
polynomials, Euler polynomials, \ldots .  
\smallskip

As useful 
as umbral calculus is in many areas of mathematics, the foundations 
of umbral calculus were not firmly established for over 
a hundred 
years. 
Vaguely speaking, the difficulty in the study is 
that such sequences do not 
observe the {\em algebra} rules of $C[x]$.  
Rota embarked on laying down
the foundation of umbral calculus during the same period of time 
as when he started the algebraic study of Baxter algebras. 
Rota's discovery is that these sequences do observe the algebra 
rules of the dual algebra (the umbral algebra), 
or in a fancier language, the coalgebra rules of $C[x]$. 
Rota's pioneer work~\cite{Ro0} was completed 
over the next decade by Rota and his 
collaborators~\cite{RKO,RR,Rom}. 
Since then, there have been a number of generalizations 
of the umbral calculus. 
\smallskip

We will give a characterization 
of umbral calculus in terms of free Baxter algebras by  
showing that the umbral algebra is the 
free Baxter algebra of weight zero on the empty set.  
We also characterize the polynomial sequences studied in 
umbral calculus in terms of operations in free Baxter 
algebras. 

%Umbral calculus and Baxter algebra are two subjects that 
%have interested Rota throughout his life time. Thus it is 
%puzzling that he did seem to aware the close connection 
%between umbral calculus and free Baxter algebras. 
%It might be due to the fact that Rota studied 
%the free Baxter algebras in the category of algebras without 
%identity. Then the umbral algebra is not a free Baxter 
%algebra; while it is in the category of algebras with 
%identity. 
%\item 

We will then use the free Baxter 
algebra formulation of the umbral calculus to give a 
generalization of the umbral calculus, called the 
{\em $\lambda$-umbral 
calculus} for each constant $\lambda$. The umbral calculus 
of Rota is the special case when $\lambda=0$. 
%
%It should be possible to give an even more generalization 
%of the umbral calculus using free Baxter algebra. 
%One can also give a formulation of the umbral calculus 
%in terms of Baxter coalgebras by combining the approach 
%in this paper and the coalgebra approach in \cite{RR,NS}. 
%We plan to carry out these projects in a future paper. 
%}

%\newpage

\subsubsection{Rota's umbral algebra}
In order to describe the binomial sequences, Rota and his 
collaborators identify $C[x]$ as the dual of the 
algebra $C[[t]]$, called the {\em umbral algebra.} 
(algebra plus the duality).
To identify $C[[t]]$ with the dual of $C[x]$, 
let $t_n=\frac{t^n}{n!},\ n\in \NN$. 
Then
\begin{equation}
t_m t_n = \binc{m+n}{m} t_{m+n},\ m,\ n\in \NN. 
\mlabel{eq:div}
\end{equation}
The $C$-algebra $C[[t]]$, together with the basis
$\{t_n\}$ is called the {\em umbral algebra}. 

We can identify $C[[t]]$ with the dual $C$-module of 
$C[x]$ by taking $\{t_n\}$ to be the dual basis of 
$\{x^n\}$. In other words, $t_k$ is defined by  
\[ t_k : C[x] \to C,\ x^n \mapsto \delta_{k,n},\ 
k,\ n\in \NN. \]
Rota and his collaborators removed the mystery of 
sequences of binomial type and Sheffer sequences by 
showing that such sequences have a simple characterization 
in terms of the umbral algebra. 

%\newpage

Let $f_n,\ n\geq 0,$ be a pseudo-basis of $C[[t]]$. 
That is, $f_n,\ n\geq 0$ are linearly independent 
and generate $C[[t]]$ as a topological $C$-module 
where the topology on $C[[t]]$ is defined by the 
filtration 
\[ F^n=\left \{\sum_{k=n}^\infty c_k t_k 
     \right \}.\]
A pseudo-basis $f_n,\ n\geq 0,$ of $C[[t]]$ 
is called a {\em divided power} pseudo-basis 
if 
\[ f_m f_n = \binc{m+n}{m} f_{m+n},\ m,\ n\geq 0. \]

\begin{thm}
\mlabel{thm:rr}
{\rm \cite{Ro0,RR} }
\been
\item
A polynomial sequence $\{p_n(x)\}$ is of binomial 
type if and only if it is the dual basis of a divided 
power pseudo-basis of $C[[t]]$.
\item
Any divided power pseudo-basis of $C[[t]]$ is of 
the form 
$ f_n(t)=\frac{f^n(t)}{n!}$ for some $f\in C[[t]]$ with 
$\ord f=1$ (that is, $f(t)=\displaystyle{\sum_{k=1}^\infty c_k t^k},\ 
c_1\neq 0$). 
\enen
\end{thm}
This theorem completely determines all polynomial 
sequences of binomial type. Algorithms to determine such sequences 
effectively have also been developed. 
See the book by Roman~\cite{Rom} for details. 
%The umbral calculus has since been extended in several 
%directions.

%\newpage

\subsubsection{$\lambda$-umbral calculus}
Our first observation is that,  with the operator 
$$P:C[[t]]\to C[[t]],\  t_n \mapsto t_{n+1},$$  
$C[[t]]$ becomes a Baxter algebra of weight zero,  
isomorphic to $\widehat{\sha}_C(C, 0)$. 
More generally,  we have 
\begin{thm} {\rm \cite{Gu3}}
\mlabel{thm:ub}
A sequence $\{f_n(t)\}$ in $C[[t]]$ is a divided power pseudo-basis 
if and only if the map 
$f_n(t)\mapsto t_n, n\geq 0$,  
defines an automorphism of 
the Baxter algebra $C[[t]]$. 
\end{thm}
\addsm

This theorem provides a link between umbral calculus and Baxter 
algebra. 
This characterization of the umbral calculus in terms of Baxter 
algebra also motivates us to study a generalization 
of binomial type sequences. 

\begin{defn}
A sequence $\{p_n (x)\mid n\in\NN \}$ of polynomials in 
$C[x]$ is a sequence of $\lambda$-binomial type
if 
\[ p_n(x+y)= \sum_{k=0}^n \lambda^k 
    \sum_{i=0}^n \bincc{n}{i}\binc{i}{k} 
    p_i(x) p_{n+k-i}(y), 
\ \forall\, y\in C,\, n\in \NN.
\]
\end{defn}
\addsm
When $\lambda=0$, we recover the sequences of binomial 
type. 
Denote 
$${e_\lambda(x)=
    \frac{e^{\lambda x}-1}{\lambda}}$$
for the series  
$${\sum_{k=1}^\infty \frac{\lambda^{k-1} x^k}{k!}}\,.$$
When $\lambda=0$, we get  
$e_\lambda(x)=x$. 
We verify that 
$${\mathfrak q}\defeq\{(e_\lambda(x))^n\}_n$$
is a sequence of $\lambda$-binomial type of $C[[x]]$. 
Let $\Cq$ be the $C$-submodule of $C[[x]]$ generated 
by elements in ${\mathfrak q}$. 

Let $f(t)$ be a power series of order 1. 
Define
$$ d_n(f)(t)= \frac{f(t)(f(t)-\lambda)
    \cdots (f(t)-(n-1)\lambda)}{n!}, \ 
    n\geq 0\,. $$
Then $P:C[[t]]\to C[[t]],\ d_n(f) \mapsto d_{n+1}(f)$ 
defines a weight $\lambda$ Baxter operator on $C[[t]]$. 
Such a sequence is called a {\em Baxter pseudo-basis} 
of $C[[t]]$. 

\begin{defn}
Fix a $\lambda\in C$. 
The algebra $C[[t]]$, together with the weight $\lambda$ 
Baxter pseudo-basis 
$\{d_n(t)\}_n$, is called the {\em $\lambda$-umbral algebra}. 
\end{defn}
\addsm

As in the classical case, 
we identify $C[[t]]$ with the dual algebra of $\Cq$ by 
taking $\{d_n(t)\}$ to be the dual basis of $\{(e_\lambda(x))^n\}$. 
We then extend the classical theory of the umbral 
calculus to the $\lambda$-umbral calculus. 
In particular, 
Theorem~\ref{thm:rr} is generalized to 
\begin{thm} {\rm \cite{Gu3}}
\begin{enumerate}
\item A pseudo-basis $\{s_n(x)\}$ of $C[[x]]$ 
is of  $\lambda$-binomial type if 
and only if 
$\{s_n(x)\}$ is the dual basis of 
a Baxter pseudo-basis of $C[[t]]$. 
\item
Any Baxter pseudo-basis of $C[[t]]$ is of the 
form $\{d_n(f)\}$ 
for some $f(t)$ in $C[[t]]$ of order 1. 
\end{enumerate}
\end{thm}

\addcontentsline{toc}{section}{\numberline {}References}

%\begin{center}
%{\Large\bf REFERENCES}
%\end{center}

%\vspace{-1.75cm}

\end{document}